# Qualitative analysis of a class of SIRS infectious disease models with nonlinear infection rate

MENGQI TAN[1]

**ABSTRACT:** The existence and local stability of some non-negative equilibrium points of a class of SIRS infectious disease models with non-linear infection and treatment rates are investigated under the condition that the total population is a constant. The qualitative theory of differential equations was used to demonstrate that the endemic equilibrium point of the system is either a stable equilibrium, an unstable equilibrium or a degenerate equilibrium under different circumstances. Subsequently, the local stability of the non-negative equilibrium point of the system is analyzed. Finally, the bifurcation theory is used to prove that the system takes the natural recovery growth rate as the parameter of the saddle-node branching, and the conditions for the existence of the model saddle-node branching are given.

## 0 Introduction

Infectious diseases are diseases that are transmitted between humans or animals by pathogens that infect humans or other organisms, and are usually characterized by high contagiousness. The danger of transmissible diseases has prompted people to study their transmission patterns and control them, so as to establish a large number of mathematical models and conduct in-depth research on the dynamics of the models. In 1927, Kermack et al. proposed the SIR bin model, and according to different bins, different models were established, including SIR, SIS, SIRS and SEIRS models. In these models, infection rate is a very important influence factor. For the infection rate, Capasso and Serio introduced a saturated incidence function $\lambda SI/(1+\alpha I)$ in the literature[1], Liu et al. studied a more general incidence function. Driessche[2] studied the nonlinear incidence $\beta I(1+vI^{k-1})S$. In the literature[2], the SIS model was studied and it was found that backward branching phenomenon can occur in this model and local and global stabilization of equilibrium points were obtained. In addition to the incidence rate, treatment is also one of the important factors affecting the outbreak of disease. Wang and Ruan[4] investigated the SIR model with a constant treatment rate $T(I)=r$, which can occur in saddle-node branching, Hopf branching, and homozygous branching. The model was used as a basis for the development of the SIR model, and it has been shown to be a useful tool for the development of the SIR model. Subsequently, Wang proposed a new treatment function $T(I)=\begin{cases} rI, 0 < I \leq I_0 \\ k, I > I_0 \end{cases}$ in[3]. This treatment function implies that the treatment rate is proportional to the number of infected individuals when and only when the treatment capacity has not been maximized, otherwise the researcher will take the maximum treatment capacity.

Based on the results of previous research, this paper proposes a class of SIRS models with nonlinear infection rates and treatment functions.

$$\begin{cases} \dfrac{dS}{dt} = A - dS - \lambda I(1+vI)S + \theta R, \\ \dfrac{dI}{dt} = \lambda I(1+vI)S - (d+\mu)I - T(I), \\ \dfrac{dR}{dt} = \mu I - dR + T(I) - \theta R. \end{cases} \quad (0.1)$$

where S is the susceptible population, I is the infected population, R is the recovered population, $A$ is the birth rate, $d$ is the mortality rate, $\lambda$ is a constant of proportionality, $\mu$ is the natural recovery rate, $\varepsilon$ is the disease lethality rate, $\theta$ is the immunization loss rate, and $T(I)$ is a treatment function. $\lambda I(1+vI)S$ refers to nonlinear incidence.

# 1 Model building

Based on the classical Kermack-Mckendric [9] model, the following model was developed, taking into account the incidence of nonlinearity and the effect of the treatment function on the disease:

$$\begin{cases} \dfrac{dS}{dt} = A - dS - \lambda I(1+vI)S + \theta R, \\ \dfrac{dI}{dt} = \lambda I(1+vI)S - (d+\mu)I - T(I), \\ \dfrac{dR}{dt} = \mu I - dR + T(I) - \theta R. \end{cases} \quad (1.1)$$

To simplify the model, the three equations of the system are added together to have $d(S+I+R)/dt = A - d(S+I+R)$, where $S+I+R = N$, thus we have $N(t) \to A/d$, when $dN/dt = A - dN$, $t \to +\infty$, so we can discuss the system consisting of only the last two equations:

$$\begin{cases} \dfrac{dI}{dt} = \lambda I(1+vI)S - (d+\mu)I - T(I), \\ \dfrac{dR}{dt} = \mu I - dR + T(I) - \theta R. \end{cases} \quad (1.2)$$

# 2 Dynamical form of the model equilibrium point

In this section, the dynamical form of the system will be studied. In the system (1.2), , X, Y denote the densities of the prey and predator populations, respectively, and are not defined on $x = 0$. From (1.2), the positive half-axis of the $x$ axis is the invariant set of the system. Considering the biological significance, the morphology of the system is discussed in the first quadrant and the positive half-axis of the $x$ axis in this paper.

## 2.1 Existence of equilibrium points

When $0 < I \leq I_0$, the system is:

$$\begin{cases} \dfrac{dI}{dt} = \lambda I(1+vI)S - (d+\mu)I - rI, \\ \dfrac{dR}{dt} = \mu I - dR + rI - \theta R. \end{cases} \quad (2.1)$$

For ease of computation, let $k = d + \theta$, $I = kx/\lambda$, $R = ky/\lambda$, $t = \tau/k$. For notational convenience, we still write $\tau$ as $t$, and the system becomes:

$$\begin{cases} \dfrac{dx}{dt} = x(1+mx)(B-x-y) - ex, \\ \dfrac{dy}{dt} = qx - y. \end{cases} \quad (2.2)$$

where $m = vk/\lambda$, $B = \lambda A/dk$, $e = (d+\mu+r)/k$, $q = (\mu+r)/k$. The equilibrium point of the system is the solution of the following equation:

$$\begin{cases} x(1+mx)(B-x-y)-ex=0, \\ qx-y=0. \end{cases} \quad (2.3)$$

Solving the system of equations (2.3) yields the following theorem:

**Theorem 1:** For any $m>0$, $B>0$, we have

(1) The system (1.3) has no positive equilibrium point when $R_0<1$, $m\leq \dfrac{1+q}{B}$.

(2) The system (1.3) has two positive equilibrium points $E_1(x_1,y_1)$ and $E_2(x_2,y_2)$ when $1-\dfrac{(1+q-Bm)^2}{4me(1+q)}<R_0<1$, $m>\dfrac{1+q}{B}$.

(3) The system (1.3) has a dual positive equilibrium point $E^*(x^*,y^*)$ when $1-\dfrac{(1+q-Bm)^2}{4me(1+q)}=R_0<1$, $m>\dfrac{1+q}{B}$.

(4) The system (1.3) has no positive equilibrium points when $R_0<1-\dfrac{(1+q-Bm)^2}{4me(1+q)}$, $m>\dfrac{1+q}{B}$.

(5) The system (1.3) has a positive equilibrium points when $R_0\geq 1$.

**Proof:** From the second equation $y=qx$, substituting into the first equation we have:
$$x[(1+mx)(B-x-qx)-e]=0,$$

That is, there exists a disease-free equilibrium point $E(0,0)$, and the local disease equilibrium point is determined by the following quadratic equation:
$$h(x)=m(1+q)x^2+(1+q-Bm)x+e-B=0,$$

where $h(0)=B-e$, and the axes of symmetry is $x=(Bm-q-1)/(2m(1+q))$, $\Delta=(1+q-Bm)^2-4m(1+q)(e-B)$, and the two zeros are $x_i=\dfrac{(Bm-q-1\pm\sqrt{\Delta})}{2m(1+q)}$,

**Case 1:** Considering the case $B=e$, we have $h(0)=0$, and it is clear that at this point $\Delta>0$. If $q<1-Bm$, then $Bm-q-1>0$, at this point the axis of symmetry is positive, and the equation has a unique zero $E_2(x_2,y_2)$. If there is $q\geq 1-Bm$, that is $Bm-q-1\leq 0$, which implies the axis of symmetry is negative at this point, and the equation has a unique negative zero, i.e., there is no positive equilibrium point for the system.

**Case 2:** Considering the case $B>e$, we have $h(0)<0$, and it is clear that at this point $\Delta>0$, so the equation has only one positive zero, i.e., the system has an equilibrium point $E_2(x_2,y_2)$.

**Case 3:** Consider the case $B<e$, we have $h(0)>0$. If $q<1-Bm$, , then $Bm-q-1>0$, and the axis of symmetry is positive. Consider the case where $e\leq B+\dfrac{(1+q-Bm)^2}{4m(1+q)}$, we have $\Delta\geq 0$, where there are two positive zeros of the equation. That is, the system has two positive equilibrium points $E_1(x_1,y_1)$. and $E_2(x_2,y_2)$. When $e>B+\dfrac{(1+q-Bm)^2}{4m(1+q)}$, the equation has a positive zero. That is, the system has a positive equilibrium point $E^*(x^*,y^*)$. If $q\geq 1-Bm$, then $Bm-q-1\leq 0$, .,then the axis of symmetry is negative, and there are no positive zeros in the equation. Let $y_1=qx_1$, $y_2=qx_2$, $y^*=qx^*$.

If $0 < x \leq x_0$, then $E_i$ is the endemic equilibrium of the system (1.2). So let's first consider the case $x_1 > x_0$ then we need to satisfy $(Bm - q - 1 - \sqrt{\Delta})/2m(1+q) > x_0$ to get $0 < Bm - q - 1 - 2m(1+q)x_0 > \sqrt{\Delta}$. From $0 < Bm - q - 1 - 2m(1+q)x_0$, it follows that

$$R_0 > (q+1)/me + 2x_0(1+q)/e \doteq \hat{P}_1. \tag{2.4}$$

Similarly, squaring $Bm - q - 1 - 2m(1+q)x_0 > \sqrt{\Delta}$ yields

$$R_0 < \frac{1}{1+mx_0} + \frac{m(1+q)}{e(1+mx_0)}x_0^2 + \frac{q+1}{1+mx_0}x_0 \doteq \hat{P}_2. \tag{2.5}$$

Therefore, $x_1 < x_0$ holds if and only if (2.4) and (2.5) are both true. It follows that $x_1 \leq x_0$ holds when $R_0 \leq \hat{P}_1$ or $R_0 \geq \hat{P}_2$. Similarly analyze the case $x_2 < x_0$. It can be shown that $x_2 > x_0$ holds when (2.4) is satisfied or $\hat{P}_2 < R_0 \leq \hat{P}_1$. That is, when $R_0 \leq \min\{\hat{P}_1, \hat{P}_2\}$, $x_2 \leq x_0$. If $\hat{P}_1 < \hat{P}_2$, i.e., it is necessary to satisfy

$$(q+1)/me + 2x_0(1+q)/e < \frac{1}{1+mx_0} + \frac{m(1+q)}{e(1+mx_0)}x_0^2 + \frac{q+1}{1+mx_0}x_0.$$

Simplifying this equation gives $me - q - 1 > m^2(3+q)x_0^2 + m(3q + 3 - qe - q)x_0$. For $x^* < x_0$, that is $Bm - q - 1 > 2m(1+q)x_0$. If $(1+q)/B > (2mqx_0 + q + 1)/(B - 2x_0)$, the simplification gives $-2(1+q)x_0 > 2Bmx_0$. This is obviously not true, so $(1+q)/B > (2mqx_0 + q + 1)/(B - 2x_0)$ is always true.

From this we obtain Theorem 2:

**Theorem 2:** In the case of Theorem 1, consider the interrupted point x, with the following result:

(1) If $me - q - 1 > m^2(3+q)x_0^2 + m(3q + 3 - qe - q)x_0$, $R_0 \leq \hat{P}_1$, then both $E_1$ and $E_2$ exist.

(2) If $me - q - 1 > m^2(3+q)x_0^2 + m(3q + 3 - qe - q)x_0$, $\hat{P}_1 < R_0 < \hat{P}_2$ then both $E_1$ and $E_2$ exist; $R_0 \geq \hat{P}_2$, then $E_1$ exists and $E_2$ does not.

(3) If $me - q - 1 \leq m^2(3+q)x_0^2 + m(3q + 3 - qe - q)x_0$, then $E_1$ exists. If $R_0 \geq \hat{P}_2$, then $E_2$ exists. If $R_0 < \hat{P}_2$, then $E_2$ does not exist.

Consider the case $I > I_0$, at which point the system (1.2) becomes:

$$\begin{cases} \dfrac{dI}{dt} = \lambda I(1+vI)S - (d+\mu)I - n, \\ \dfrac{dR}{dt} = \mu I - dR + n - \theta R. \end{cases} \tag{2.6}$$

After variable substitution, the system becomes:

$$\begin{cases} \dfrac{dx}{dt} = x(1+mx)(B - x - y) - gx - f, \\ \dfrac{dy}{dt} = px - y + f. \end{cases} \tag{2.7}$$

where $m = \dfrac{vk}{\lambda}$, $B = \dfrac{\lambda A}{dk}$, $g = \dfrac{d+\mu}{k}$, $p = \dfrac{\mu}{k}$, $f = \dfrac{n\lambda}{k^2}$. The equilibrium point of the system

(2.7) is determined by the following equation:
$$\begin{cases} x(1+mx)(B-x-y) - gx - f = 0, \\ px - y + f = 0. \end{cases} \quad (2.8)$$

Solving the system of equations (2.8) leads to the following conclusion:

**Theorem 3:**

(1) If $p = Bm - mf - 1$, $B > g + f$, $-27m(1+p)f^2 = 4(g+f-B)^3$, then the system (2.8) has a dual equilibrium point $E_3(x_3, y_3)$; In particular, if $f > m(1+p)x_0^3$, then $E_3(x_3, y_3)$ is the endemic equilibrium of system (1.2).

(2) If $p = Bm - mf - 1$, $B > g + f$, $-27m(1+p)f^2 > 4(g+f-B)^3$ the system (2.8) has two positive equilibrium points $E_4(x_4, y_4)$ and $E_5(x_5, y_5)$.

**Proof:** From the second equation in (2.8), we know that $y = qx$, and substituting into the first equation, we have $m(1+p)x^3 + (1+p+mf-Bm)x^2 + (g+f-B)x + f = 0$.

It is known that the endemic equilibrium point is determined by the following one-dimensional cubic equation:
$$F(x) = m(1+p)x^3 + (1+p+mf-Bm)x^2 + (g+f-B)x + f = 0,$$

Then $\dfrac{1}{m(1+p)} F(x)$ is $x^3 + \dfrac{(1+p+mf-Bm)}{m(1+p)} x^2 + \dfrac{g+f-B}{m(1+p)} x + \dfrac{f}{m(1+p)} = 0.$

Due to the complexity of solving one-dimensional cubic equations, here we present only the more special cases. If $p = Bm - mf - 1$, i.e. $1+p+mf-Bm = 0$, then the equation is $x^3 + sx + t = 0$. Where $s = (g+f-B)/m(1+p)$, $t = f/m(1+p) > 0$.

From literature[11], the equation exists with two equal positive real roots if $t^2/4 + s^3/27 = 0$ i.e. The equation holds if and only if $g + f < B$. And from $1 + p + mf - Bm = 0$, we know that $1 + p + mf - Bm = 0$. Thus we have $-27m^2 f^2 (B-f) = 4(g+f-B)^3$. The root at this point is $x_3 = -\sqrt[3]{-t/2}$. If $t^2/4 + s^3/27 > 0$, then the equation has two unequal positive real roots.

If $x > x_0$, then $E_i (i = 3, \cdots, 5)$ is an endemic equilibrium point of the system (1.2). Due to the complexity of the roots of a cubic equation, we will only discuss the case of dual roots here.

If $x_3 > x_0$, then $-\sqrt[3]{-\dfrac{t}{2}} > 0.$, i.e., $\sqrt[3]{\dfrac{f}{m(1+p)}} > x_0 > 0$, and thus $f > m(1+p)x_0^3$. That is, under our assumptions, if $f > m(1+p)x_0^3$, then $E_3(x_3, y_3)$ is an endemic equilibrium point of the system (1.2).

## 2.2 Dynamical states at equilibrium points

The dynamical state at the equilibrium point is the key to further study of the model, and analyzing the stability of the system at the equilibrium point can suggest solutions to the practical problems corresponding to the model. For this purpose, the eigenvalues of the Jacobi matrix of the system at the equilibrium point are considered.

Note: This section considers the stability of the equilibrium point when the treatment rate has not yet reached the maximum treatment capacity.

The Jacobi matrix of the system (1.3) at the equilibrium point $E(x, y)$ is
$$J = \begin{pmatrix} B - 2x - qx + 2Bmx - 3mx^2 - 2mqx^2 - e & -x - mx^2 \\ q & -1 \end{pmatrix},$$

The characteristic equation is given by $\lambda^2 - tr(J)\lambda + \det(J) = 0$. $tr(J)$ and $\det(J)$ denote the trace and determinant of the matrix. The expressions are as follows:

$$tr(J) = B - 2x - qx + 2Bmx - 3mx^2 - 2mqx^2 - e - 1.$$
$$\det(J) = 3m(1+q)x^2 + 2(1+q-Bm)x + e - B.$$

Note 1: The $\det(J|_E)$ symbol is determined by $p(x) = x(2m(1+q)x + 1 + q - Bm)$.

Note 2: The $tr(J|_E)$ symbol is determined by $q(x) = ((q+1)^2 - Bm)x + (e-B)(q+2) - (q+1)$.

### 2.2.1 Stability of disease-free equilibrium points

The Jacobi matrix of the system at the disease-free equilibrium point is

$$J|_{E_0} = \begin{pmatrix} B-e & 0 \\ q & -1 \end{pmatrix},$$

The eigenvalues of this matrix are $\lambda_1 = -1 < 0$, $\lambda_2 = B - e$. The determinant and trace of theJacobi matrix $J|_{E_0}$ are respectively: $\det(J|_{E_0}) = -(B-e)$, $Tr(J|_{E_0}) = B - e - 1$. By the sign of $\det(J|_{E_0})$ the following theorem follows.

**Theorem 4:**

(1) If $B > e$, the disease-free equilibrium point $E_0 = (0,0)$ is unstable; if $B < e$, the disease-free equilibrium point $E_0 = (0,0)$ is locally and progressively stable.

(2) If $B = e$, the disease-free equilibrium point $E_0$ is a saddle point.

**Proof:** (1) From reference[5], if $B > e$, the disease-free equilibrium point $E_0 = (0,0)$ is unstable; if $B < e$, the disease-free equilibrium point $E_0 = (0,0)$ is locally asymptotically stable.

(2) When $B = e$, the two eigenvalues of the system are 0,-1. A Taylor expansion of this system at the origin yields the following system:

$$\begin{cases} \dfrac{dx}{dt} = a_{10}x + a_{01}y + a_{20}x^2 + a_{11}xy + a_{02}y^2 + o((u,v)^3), \\ \dfrac{dy}{dt} = b_{10}x + b_{01}y + b_{20}x^2 + b_{11}xy + b_{02}y^2 + o((u,v)^3). \end{cases} \quad (2.9)$$

Where $a_{10} = 0$, $a_{01} = 0$, $a_{11} = -1$, $a_{20} = -2 + 2Bm$, $a_{02} = 0$, $b_{10} = q$, $b_{01} = -1$, $b_{11} = 0$, $b_{20} = 0$, $b_{02} = 0$.

Then the following system is obtained:

$$\begin{cases} \dfrac{dx}{dt} = 2(Bm-1)x^2 - xy + o((u,v)^3), \\ \dfrac{dy}{dt} = qx - y. \end{cases} \quad (2.10)$$

From reference[6], if $Bm - 1 < 0$, the point is a stable saddle point; if $Bm - 1 > 0$, the point is an unstable saddle point.

**Theorem 5:** The disease-free equilibrium of the system e is globally asymptotically stable if the following conditions hold.

(1) $e - B > 0$, $m \leq (1+q)/B$.

(2) $B < 1 < g, mB < 1$.

**Proof:** Under condition (1), the system does not have endemic equilibrium points by Theorem 1, while the disease-free equilibrium point e is locally asymptotically stable by Theorem 4, so the disease-free equilibrium point e is globally asymptotically stable.

In condition (2), let $P(x,y) = Bx - x^2 - xy + Bmx^2 - mx^3 - mx^2y - ex$, $Q(x,y) = qx - $

$y$. When $I \leq I_0$, take the Dulac function $D(x,y) = 1/xy, x \leq x_0 = \lambda I_0/k$. we have

$DP = [(1+mx)(B-x-y) - e/y]/y$, $DQ = q/y - 1/x$, $\dfrac{\partial(DP)}{\partial x} + \dfrac{\partial(DQ)}{\partial y} = \dfrac{1}{y^2}(-my^2$

$+(mB-1-2mx)y-q)$, Since $mB < 1$, it is constant that $\dfrac{\partial(DP)}{\partial x} + \dfrac{\partial(DQ)}{\partial y}$ is less than 0.

Similarly, when $I > I_0$, take the Dulac function $D(x,y) = 1/x_0 y, x > \lambda I_0/k = x_0$, we have

$DP = [x(1+mx)(B-x-y) - gx - f]/x_0 y$, $DQ = (px - y + f)/x_0 y$.

$\dfrac{\partial(DP)}{\partial x} + \dfrac{\partial(DQ)}{\partial y} = \dfrac{1}{x_0 y^2}(-y^2 + (-3mx^2 + 2m(B-1)x - x + (B-g)]y - px - f)$.

Also by $B < 1 < g$, thus the above equation is less than 0 constant.

From the literature[10], the system does not have a limit loop. Thus in summary $E_0$ is globally asymptotically stable.

### 2.2.2 Endemic disease equilibrium stabilization

In this subsection, we discuss the stability of the endemic equilibrium point, and due to space constraints, here we only discuss the case of stability at the dichotomous equilibrium point $E^*$. Let $\varepsilon = q(x^*) = ((q+1)^2 - Bm)x^* + (e-B)(q+2) - (q+1)$.

**Theorem 6:** When $\varepsilon \neq 0$, $E^*$ is a saddle node.

**Proof:** Substituting $x^* = (Bm - q - 1)/2m(1+q)$ for $p(x^*)$, there is: $p(x) = x^*(2m(1+q)x^* + 1 + q - Bm) = 0$. Therefore $\det(J|_{E^*}) = 0$, i.e. this equilibrium is degenerate. If $\varepsilon \neq 0$, $q(x^*) \neq 0$, i.e. $tr(J|_{E^*}) \neq 0$, then the system has only one zero eigenvalue at this equilibrium point, by $u = x - x^*, v = y - y^*$. Shifting the equilibrium to the origin, the following system is obtained:

$$\begin{cases} \dfrac{du}{dt} = (u+x^*)(1+m(u+x^*))(B-u-x^*-v-y^*) - e(u+x^*), \\ \dfrac{dv}{dt} = qu - v. \end{cases} \quad (2.11)$$

A Taylor expansion of the above system at the origin is obtained:

$$\begin{cases} \dfrac{du}{dt} = a_{10}u + a_{01}v + a_{20}u^2 + a_{11}uv + a_{02}v^2 + o((u,v)^3), \\ \dfrac{dv}{dt} = b_{10}u + b_{01}v + o((u,v)^3). \end{cases} \quad (2.12)$$

Where $a_{10} = B - 2x^* - y^* + 2Bmx^* - 3mx^{*2} - 2mx^*y^* - e$, $a_{01} = -x^* - mx^{*2}$, $a_{20} = -2 + 2Bm - 6mx^* - 2my^*$, $a_{02} = 0$, $a_{11} = -1 - 2mx^*$, $b_{10} = q$, $b_{11} = -1$.

The linearization matrix of the system at $E^*$ is:

$$J|_{E^*} = \begin{pmatrix} a_{10} & a_{01} \\ q & -1 \end{pmatrix},$$

and $-a_{10} - qa_{01} = 0$, $a_{10} - 1 \neq 0$.

Doing a non-singular linear transformation $u = u_1 - a_{01}v_1$, $v = qu_1 + v_1$, and a time transformation $dt = -(a_{10} + 1)d\tau$ for the system (2.12), the following system is obtained:

$$\begin{cases} \dfrac{du_1}{dt} = c_{20}u_1^2 + o((u_1,v_1)^3), \\ \dfrac{dv_1}{dt} = d_{01}v + o((u,v)^2). \end{cases} \tag{2.13}$$

where $c_{20} = \dfrac{1}{1-a_{10}}(a_{20} + qa_{11})$, $d_{01} = \dfrac{b_{10}}{1-a_{10}} \neq 0$.

From the literature[7], the equilibrium point is a saddle node. When $a_{10} > 1, c_{20} < 0$, and the time transformation is negative, so the equilibrium point $E^*$ is the saddle node of repulsion; when $a_{10} < 1$, $c_{20} > 0$, and the time transformation is negative, so the equilibrium point $E^*$ is the saddle node of attraction.

## 3 Branching Analysis

Studying the kinetic behaviors such as branching of the SIRS model with nonlinear incidence is beneficial to human beings for better prevention of infectious diseases. In this subsection, we will determine the possible existence of branches of the system (1.2), including saddle-node branch and hopf branch, based on the analysis results of the dynamical state of the equilibrium point in the case of the previous subsection. Here we only give the conditions for the existence of saddle-node branches.

### 3.1 Saddle-node bifurcation

It follows from Theorem 2 that if $e > B$, $m > \dfrac{1+q}{B}$, $e < B + \dfrac{(1+q-Bm)^2}{4m(1+q)}$, the system has two internal equilibria $E_1$ and $E_2$; if $e > B$, $m > \dfrac{1+q}{B}$, $e > B + \dfrac{(1+q-Bm)^2}{4m(1+q)}$, the system has no internal equilibrium point; if $e > B$, $m > \dfrac{1+q}{B}$, $e = B + \dfrac{(1+q-Bm)^2}{4m(1+q)}$, there is a unique positive equilibrium point $E^*$.

So let $e_{SN} = B + (1+q-Bm)^2/4m(1+q)$, when $e$ varies around $e_{SN}$, the number of boundary equilibrium points changes, giving rise to saddle-node branches.

**Theorem 9:** When $\varepsilon \neq 0$, the system (1.2) undergoes a saddle-node bifurcation at $E^*$ with a critical bifurcation parameter $e = e_{SN} = B + (1+q-Bm)^2/4m(1+q)$.

**Proof:** The Jacobian matrix of the system (1.2) at $E^* = (x^*, y^*)$ is:

$$J|_{E^*} = \begin{pmatrix} B - 2x^* - qx^* + 2Bmx^* - 3mx^{*2} - 2mqx^{*2} - e & -x^* - mx^{*2} \\ q & -1 \end{pmatrix},$$

Because $\varepsilon \neq 0$, so $J(E^*)$ has a zero eigenvalue, let $v$, $w$ be the eigenvectors corresponding to the zero eigenvalues of $J(E^*)$ and $J^T(E^*)$, respectively. Then we can obtain:

$$v = \begin{pmatrix} v_1 \\ v_2 \end{pmatrix} = \begin{pmatrix} 1 \\ q \end{pmatrix},$$

Let
$$w = \begin{pmatrix} w_1 \\ w_2 \end{pmatrix} = \begin{pmatrix} q \\ x^* + mx^{*2} \end{pmatrix},$$

$$F(x,y) = \begin{pmatrix} \dot{x} \\ \dot{y} \end{pmatrix} = \begin{pmatrix} F_1 \\ F_2 \end{pmatrix} = \begin{pmatrix} x(1+mx)(B-x-y)-ex \\ qx-y \end{pmatrix},$$

We obtain:
$$F_e(E^*, e_{SN}) = \begin{pmatrix} -x^* \\ 0 \end{pmatrix},$$

$$D^2 F_e(E^*, e_{SN})(v,v) = \begin{pmatrix} \frac{\partial^2 F_1}{\partial x^2} v_1^2 + 2\frac{\partial^2 F_1}{\partial x \partial y} v_1 v_2 + \frac{\partial^2 F_1}{\partial y^2} v_2^2 \\ \frac{\partial^2 F_2}{\partial x^2} v_1^2 + 2\frac{\partial^2 F_2}{\partial x \partial y} v_1 v_2 + \frac{\partial^2 F_2}{\partial y^2} v_2^2 \end{pmatrix} = \begin{pmatrix} \frac{2Bm - 3q - 2 - q^2}{1+q} \\ 0 \end{pmatrix}.$$

So we get:
$$w^T F_e(E^*, e_{SN}) = -qx^* \neq 0,$$

$$w^T [D^2 F_e(E^*, e_{SN})(v,v)] = \frac{2Bmq - 3q^2 - 2q - q^3}{1+q} \neq 0.$$

By Sotomayor's theorem[8], it is known that the system (1.2) has a saddle-node bifurcation at the equilibrium point $E^*$.

## 4 Conclusion

In this paper, we study a class of SIRS models with nonlinear incidence and treatment functions. First, the existence case of the internal equilibrium point of the system (1.2) is analyzed, and it is obtained that, under certain conditions, the system has one or two internal equilibrium points. Secondly, the stability of the disease-free equilibrium point is analyzed and it is obtained that the equilibrium point is a saddle point. Subsequently, after theoretical analysis, it is obtained that the dual positive equilibrium point of the system can be stable, unstable or degenerate. Finally, using the bifurcation theory, it is obtained that the system undergoes saddle-node branching at the dual equilibrium point when the parameter varies in the neighborhood.

[1]College of Mathematics and Statistics, Chongqing University, Chongqing, 401331, P. R. China
   Email address: tmqhhy@163.com